\title{On Blocking Sets of Affine Spaces}

\author{Ara Aleksanyan\\
Department of Applied Mathematics,\\  
Yerevan State University, Yerevan 375049, Armenia\\
and\\
Mihran Papikian \thanks{e-mail: papikian@math.lsa.umich.edu}\\
Department of Mathematics,\\ 
University of Michigan, Ann Arbor, Michigan 48109\\}

\date{} %Prevent the date from being printed
\documentstyle [11pt]{article}

\newtheorem{theorem}{Theorem}[section]
\newtheorem{proposition}[theorem]{Proposition}

\newtheorem{lemma}[theorem]{Lemma}

\newcommand{\ignore}[1]{}

\begin{document}
\maketitle
\setcounter{page}{1}
\begin{abstract}

\noindent
Let $B(k, n, q)$ be a set of points of an affine space $AG(n, q)$ such that 
$C \cap B(k, n, q) \neq \emptyset$ for any affine subspace $C$ of dimension k 
in $AG(n, q)$. Bounds on the possible minimal order of  $B(1, n, q)$ and 
$B(2, n, 2)$ are obtained. Several examples are discussed.
\end{abstract}

\section{Introduction}

We consider the affine space $AG(n, q)$ over finite field $F_q$, $q$ is a power of a prime, 
and consequently can associate
$AG(n, q) \cong F_q^n$ where $F_q^n$ is a $n$-dimensional 
linear space over the field $F_q$, and $F_q^n$
can be uniquely associated with $F_{q^n}$ - the field extension of $F_q$ of degree $n$. 

Given an affine space $AG(n, q)$ one is interested in finding the minimal cardinality 
of a set $S$ of points, such that each $k$-dimensional affine subspaces in $AG(n, q)$ 
contains at least one point of $S$.

In this paper we consider the case $k=1$ in general, and the case $k=2$ when the affine 
space is over $F_2$.

We should make a remark that when the affine space is replaced by a projective space 
$PG(n, q)$ over $F_q$ then the similar problem of finding the minimal cardinality of a set of 
points having nonempty intersection with every $k$-dimensional subplane has a rather 
straightforward solution. R. Bose and R. Burton proved the following result in 
\cite{Bose}:
\begin{theorem}[R. Bose and R. Burton, 1966]
If S is a set of points in PG(n, q) which has a non-empty intersection with every k-dimensional 
subplane, then the number of points in S is greater than or equal to $(q^{n+1-k}-1)/(q-1)$. Equality 
holds if, and only if, S is a (n-k)-dimensional subplane.
\end{theorem}

Note that this result not only gives the minimal possible order for $S$, but also uniquely characterizes 
such sets as subplanes of certain dimension.

The situation seems to be much more complicated for affine spaces. Complete answer is known only for the 
case of point sets intersecting all affine subspaces of codimension 1. Jamison \cite{Jamison}, 
and independently Brouwer and Schrijver \cite{Brouwer}, have shown that the order of a point set 
intersecting all hyperplanes (affine subspaces of codimension 1) is not less than $n(q-1) + 1$, and 
that the bound is sharp. No classification result is known.

We define a {\em blocking set} $S$ of a finite affine space $AG(n, q)$ to be a point set in $AG(n, q)$ 
which meets every line. A blocking set which does not contain a smaller blocking set is called 
{\em irreducible}. 
Some variations of this definition are possible, for example, one may require that $S$ doesn't contain 
a line in $AG(n, q)$ (cf. Mazzocca and Tallini \cite{Mazzocca}). 
When $q=2$ the blocking sets are degenerate: since any two points in $AG(n, 2)$ form a line, to intersect 
all lines, $S$ must have at least $2^n -1$ points. So for $AG(n, 2)$ we, instead, will consider the point sets 
which have nonempty intersection with every affine subspace of dimension 2, and 
call them {\em 2-blocking sets}.
 
A blocking set in $AG(n, q)$ (or 2-blocking set for $AG(n, 2)$) of the least possible order 
will be called {\em minimal}. 
Minimal blocking sets are automatically irreducible, but not conversely as will be shown later. 

In our considerations it is more convenient to analyze the complement of a blocking set in $AG(n, q)$. 
We will call $S=AG(n, q)\backslash S'$ {\em linefree} ({\em 2-free} for $AG(n, 2)$), 
where $S'$ is a blocking set (2-blocking set). The complement of an 
irreducible blocking set will be called {\em complete}, and the complement of a minimal blocking set will 
be called {\em maximal} linefree set in $AG(n, q)$. 
\vspace{0.4cm}

The main purpose of this paper is to obtain bounds on the cardinality of minimal 2-blocking sets 
in $AG(n, 2)$, and improve the bounds 
$q^{n-1} + q^{n-2} + \cdots + 1 \leq |S| \leq q^n - (q-1)^n $ for a minimal blocking set $S$ in 
$AG(n, q)$ when $q=3$. This is done in Sections 3 and 4. 
In Section 2 we give a general construction for blocking sets from which the bounds 
$q^{n-1} + q^{n-2} + \cdots + 1 \leq |S| \leq q^n - (q-1)^n $ follow. 
\vspace{0.4cm}

For Sections 3 and 4 we need the following useful characterization of affine subspaces in $AG(n, q)$.

\begin{proposition}
\label{affine sum}
The set ${\cal C}=\{c_1, c_2, \dots, c_m \}$ is an affine subspace of $F_q^n$, that is 
${\cal C}= \alpha + {\cal L}$ for some linear subspace ${\cal L}$ and vector $\alpha$ of $F_q^n$,
if, and only if, $\sum\limits_{i=1}^m\lambda _i c_i \in {\cal C}$ for every 
$\lambda _1,...,\lambda _m \in F_q$ such that $\sum\limits_{i=1}^m\lambda _i=1$.
\end{proposition}

\section{General construction of linefree sets}

The next proposition can be interpreted as the analog for affine spaces of Bose-Burton result. 

\begin{proposition}
The complement of the union of n hyperplanes in general position in $AG(n, q)$ forms a complete 
linefree set. Such set $S$ contains $(q-1)^n$ points.
\end{proposition}

\noindent
{\tt Proof}

After an appropriate nondegenerate affine transformation we may assume that each hyperplane $\Pi_k$ in 
the proposition is $span(e_1,\dots, e_{k-1},$ $e_{k+1},\dots,e_n)$, where $e_i$'s form the orthonormal 
basis. Then the complement in $AG(n, q)$ of $\bigcup_{k=1}^n \Pi_k$ is the set $S$ of vectors in 
$AG(n, q)$ without zero coordinates.

In $AG(n, q)$ a line $\ell$ through the point $A=(a_1, \dots, a_n)$ in direction 
$B=(b_1, \dots, b_n) \neq (0,0, \dots,0)$ is the set of points 
$\ell=\{ (a_1 + \lambda b_1, a_2 + \lambda b_2, \dots , a_n + \lambda b_n) | \lambda \in F_q \}$. 
Suppose $b_i \neq 0$. Then there exists $\lambda$ such that $a_i + \lambda b_i = 0$. So each line 
in $AG(n, q)$ has a point with zero coordinates, and we get that $S$ is linefree.

To prove the completeness of $S$ suppose $A=(a_1, \dots, a_n)$ is a point not in $S$, hence $a_i=0$ for 
at least one $a_i$. Put
$$
b_i = 
\left\{
\begin{array}{cl}
1 & \mbox{if $a_i=0$}\\
0 & \mbox{if $a_i \neq 0$}
\end{array}
\right.
$$
Then 
$$
a_i + \lambda b_i =
\left\{
\begin{array}{cl}
\lambda & \mbox{if $a_i=0$}\\
a_i & \mbox{if $a_i \neq 0$}
\end{array}
\right.
$$
Easy to see that the line through $A$ in direction $B=(b_1, \dots, b_n)$ contains only one point 
(namely $A$ itself) with one or more zero coordinates. Hence $S$ is complete. 

That $S$ contains $(q-1)^n$ points is trivial. $\hfill{\Box}$
\vspace{0.4cm}

Now consider an arbitrary linefree set $S$. 
Let $P \in S$ then any line through $P$ has at least one point not in $S$. There are $\frac{q^n-1}{q-1}$ 
such lines, so $S$ cannot contain more than $q^n - \frac{q^n-1}{q-1}$ points. Summarizing, if $S$ is a 
maximal linefree set in $AG(n, q)$ then

\begin{equation}
\label{general bound}
(q-1)^n \leq |S| \leq q^n - \frac{q^n-1}{q-1}
\end{equation}

In the next section we obtain improvements on these bounds in case when $q=3$. Since we are able to obtain 
an improvement on the lower bound in (\ref {general bound}) the complement of the union 
of $n$ hyperplanes in general position is an example of a complete linefree set which is not maximal.

\section{Linefree sets in $AG(n, 3)$}
The aim of this section is to prove the following:
\begin{theorem}
Let $S$ be a linefree set of maximal order in $AG(n, 3)$. Write $n=4k+m, \ m \in \{ 0, 1, 2, 3 \}$, then
$$
20^k \cdot 9 \leq |S| \leq 3^{n-1}+1 \qquad m=3
$$
$$
20^k \cdot 2^m \leq |S| \leq 3^{n-1}+1 \qquad m \in \{ 0, 1, 2\}
$$
\end{theorem}

Let us start with an observation that in the case of $AG(n, 3)$ 
Proposition (\ref{affine sum}) gives an especially convenient
criteria for $S$ to be linefree:

\begin{lemma}
A set of points $S$ is linefree in AG(n, 3) if, and only if, for any distinct 
$s_1, s_2 \in S, -(s_1+s_2)\not\in S$.
\end{lemma}

If $S$ is linefree then obviously such are also $\alpha + S$ and $\alpha \cdot S$ for any 
$\alpha \in F_{q^n}\backslash\{0\}$.
So we may assume without any loss of generality that if $S$ is linefree then $0 \in S$. 
For such $S$ we must have $(s\neq 0) \in S \Rightarrow -s \not \in S$.

Let $S=\{s_1, s_2, \dots, s_m \}$ and  $-S=\{-s_1, -s_2, \dots, -s_m \}$. Then $S\cap( -S)=\{{\bf 0}\}$.
Suppose $s_1\neq 0$. Define a mapping for all $s_j \not \in \{s_1, 0\}$:

$$
\phi (s_j) = 
\left\{
\begin{array}{cl}
-(s_1 + s_j), & \mbox{when $s_1 + s_j \not \in S$}\\
-(s_i + s_j), & \mbox{when $s_1 + s_j =s_i \in S$ }
\end{array}
\right.
$$

$\phi$ is a mapping from $S\backslash\{s_1, 0\}$ to $AG(n, 3)\backslash(S\cup (-S))$. Indeed, for all
$j \neq 1$, $-(s_1 + s_j) \not \in S$; if $-(s_1 + s_j) \in -S$ then $s_1 + s_j=s_i$ for some $i$.
But then again $-(s_i + s_j) \not \in S$ and we must have $-(s_i + s_j) \not \in -S$ since
otherwise $s_i + s_j=s_k$ for some $k$ which along with $s_1 + s_j=s_i$ yields $-(s_1+s_i)=s_k \in S$ in
contrary with $S$ being linefree.

In fact, $\phi$ is an injection. There are three cases to check:

\begin{enumerate}

\item If $s_1 + s_j \not \in S$ and $s_1 + s_i \not \in S$ then, since $-(s_1 + s_i)\neq-(s_1 + s_j)$, 
$s_j$ and $s_i$ are mapped into different elements.

\item Suppose $s_j$ and $s_k$ are such that $s_1 + s_j = s_i\in S$ but $s_1 + s_k \not\in S$. Then
$s_j$ and $s_k$ cannot be mapped into the same element since otherwise $-(s_i + s_j)=-(s_1 + s_k) 
\Rightarrow s_i + s_j-s_1 = s_k$, summing with $s_1 + s_j-s_i=0$ we get $s_j=-s_k$. 
Contradiction. 

\item  Suppose $s_j$ and $s_k$ are such that $s_1 + s_j = s_i\in S$ and $s_1 + s_k = s_t\in S$. Then again
$s_j$ and $s_k$ cannot be mapped into the same element since otherwise $(s_i + s_j)=(s_t + s_k)$ which 
along with $s_j - s_k = s_i - s_t$ implies $s_j=s_k$.
\end{enumerate}

We have  $S\cap (-S)=\{{\bf 0}\}$ and $\phi(S\backslash\{0, s_1\})\cap (S\cup(-S))=\emptyset$, so
$3^n \geq |S| + (|S|-1)+(|S|-2)$, or
$$
3^{n-1}+1 \geq |S|
$$

\noindent
{\bf Example}
 
For $F_9$ this bound shows that $S=\{(0,0), (0,1),$ $(1,0), (1,1)\}$ is an example of a complete linefree
set of maximal order (compare with Section 2).

Let $F_{27}=\{\alpha^i\}_{i=0}^{25} \cup \{0\}$, where $\alpha$ is a primitive element of $F_{27}$ 
(we take $\alpha$ to be a root of $f(x)=x^3-x^2+1$ over $F_3$; see Table A in \cite{Lidl}). Then
$S=\{{\bf 0}=(0,0,0), \alpha^0=(0,0,1), \alpha=(0,1,0), \alpha^2=(1,0,0), \alpha^4=(1,2,2), 
\alpha^7=(1,0,1), \alpha^8=(1,1,2), \alpha^{18}=(0,1,1), \alpha^{23}=(2,1,2) \}$ can be shown to be a complete
linefree of maximal order (i.e. for AG(3, 3) the maximal linefree set has 9 elements). 
This means that the bound $3^{n-1}+1$ is not sharp in general. It is also worth to
note that the above constructed $S$ for $F_{27}$ has a nice algebraic structure: 
Consider the cyclomatic classes of $F_{27}$, that is the sets $\{ \beta, \beta^3, \beta^9\}, \quad \beta \in F_{27}$.
There are 11 of them, three of which are degenerate - $\{0\}, \{1\}, \{\alpha ^{13}=-1\}$. All elements
of $S$ except $\{\alpha^2, \alpha^{18}\}$ belong to different cyclomatic classes.

We tend to believe that a similar construction is possible in the general case of $F_{3^n}$, providing a lower
bound on $|S|$ close to $3^{n-1}$.  $\hfill{\Box}$
\vspace{0.4cm}

Now associate $AG(n, 3)\simeq F_{3^n}$ and consider the roots of unity of degree $3^{n/2}+1$ assuming $n=0\ (mod\ 2)$.
That is consider $\Omega\equiv \{ x \ | \ x^{3^{n/2}+1}=1\}$. $|\Omega|=3^{n/2}+1$ since 
$3^n-1=0 \ (mod\ 3^{n/2}+1) \Rightarrow x^{3^n-1}-1=0 \ (mod\ x^{3^{n/2}+1}-1) \Rightarrow$ $x^{3^{n/2}+1}-1=0$ completely
splits in $F_{3^n}$, and has no repeating roots.

We claim that $\Omega$ is linefree. Indeed, take any $x, y \in \Omega$ such that $x \neq y$, and 
suppose the contrary 
$-(x+y) \in \Omega$ (see Lemma (\ref{affine sum})). Then
$$
(-(x+y))^{3^{n/2}+1}=1,
$$
$$
(x+y)^{3^{n/2}}(x+y)=1, \quad \mbox{and since $x, y \in \Omega$}
$$
$$
(x^{-1}+y^{-1})(x+y)=1 \Rightarrow (x^{-1}y)^2+x^{-1}y+1=0 \Rightarrow
$$
$$
(x^{-1}y-1)^2=0 \Rightarrow x=y
$$

Note that $\Omega$ has fewer elements than the bound $2^n$ in Proposition \ref{general bound}. 
But being a cyclic group of 
$F_{3^n}\backslash \{0\}$ it has a better algebraic structure than the vectors with $0, 1$ coordinates
as given in the above mentioned proposition, and we can operate with $\Omega$ easier.
The following lemma is true for any finite field $F_q$ and any multiplicative subgroup $G$ of 
$F_q^*\equiv F_q\backslash \{0\}$:

\begin{lemma}
Let $g$ be a primitive element of $F_q^*$. Let $G=\{1, g^k, g^{2k}, \dots, g^{k(m-1)}\}$ be a subgroup
of $F_q^*$ of order $m$ where $m=\frac{q-1}{k}$, and let $g^{i_1}G$, $g^{i_2}G$, $g^{i_3}G$ be three 
distinct cosets of $G$.

If some $b_3 \in g^{i_3}G$ can be expressed as $b_3=-(b_1+b_2)$ where $b_1 \in g ^{i_1}G$, 
$b_2 \in g ^{i_2}G$ then any element of $g^{i_3}G$ can be expressed in that form, i.e.
$$
g^{i_3}G \subset \{-(b_1 + b_2) \ | \ b_1 \in g ^{i_1}G, \ b_2 \in g ^{i_2}G \}.
$$
\end{lemma}
 
\noindent
{\tt Proof}
Suppose $b_1=g^{t_1k+i_1} \in g^{i_1}G$, $b_2=g^{t_2k+i_2} \in g^{i_2}G$,
$b_3=g^{t_3k+i_3} \in g^{i_3}G$ and 
$$
-(g^{t_1k+i_1}+ g^{t_2k+i_2})=g^{t_3k+i_3}
$$
To get any $g^{rk+i_3}$ form
$$
-(g^{(t_1+r-t_3)k+i_1}+ g^{(t_2+r-t_3)k+i_2})=g^{rk+i_3}
$$ 
$\hfill{\Box}$
\vspace{0.4cm}

This lemma along with the observation that any coset of $\Omega$ is linefree provides a method
to form linefree sets as unions of the multiplicative cosets of $\Omega$.
\vspace{0.4cm}

Let $g$ be a primitive element of $F_{81}$ 
(we take $g$ to be a root of $f(x)=x^4-2x^3-1$ over $F_3$; see Table A in \cite{Lidl}). 
Then $\Omega=\{g^{8i}\ | \ i=0, 1, \dots , 9\}$. Using
the lemma to simplify the computations we check that 
$-(\Omega + \Omega) \equiv \{-(g^{8i}+g^{8j}) \ | \ 1 \leq i < j \leq 9\}=g^2\Omega \cup
g^5\Omega \cup g^6\Omega \cup g^7\Omega$. At once it also follows that
$-(g^4\Omega + g^4\Omega) =g\Omega \cup g^2\Omega \cup g^3\Omega \cup g^6\Omega$. 
But then it is easy to conclude that $\Omega \cup g^4\Omega$ is a linefree set with 20 elements, when the
example in Section 2 has only 16.
\vspace{0.4cm}

Suppose we have linefree sets $S_1$ and $S_2$ in $AG(n, 3)$ and $AG(m, 3)$ respectively. Then if we
form a new set $S \subset AG(n+m, 3)$ consisting of all the vectors $(a, b)$, where $a\in S_1$ 
and $b\in S_2$, $S$ will be linefree and $|S|=|S_1||S_2|$. 
Indeed, suppose $(a_1,b_1)$ and $(a_2,b_2)$ are in $S$ and
suppose $a_1\neq a_2$, then since $-(a_1+a_2) \not \in S_1 \Rightarrow -((a_1,b_1)+(a_2,b_2)) \not \in S$.

Apply this method to the linefree set of order 20 in $AG(4,3)$, order 9 in  $AG(3,3)$, 
order 4 in $AG(2,3)$ and order 2 in  $AG(1,3)$. We get that if $S$ is a linefree set of maximal
order in $AG(n, 3)$ then 
$$
|S| \geq 20^k \cdot 9 \qquad \mbox{for $n=4k+3$}
$$
$$
|S| \geq 20^k \cdot 2^m \qquad \mbox{for $n=4k+m, \quad m \in \{0,1,2\}$}. 
$$
Note that the bound in Section 2 gives only $16^k \cdot 2^m$.

\section{2-free sets in $AG(n, 2)$}

Any two points in $AG(n, 2)$ form a line , and consequently the linefree sets
are the trivial one-point sets.

Now let us investigate the first nontrivial case. If $S$ is a set of points in $AG(n, 2)$ 
not containing any affine subspace of dimension 2 (call such set {\em 2-free}) we are interested 
in finding the size of the largest such set. Using Proposition (\ref{affine sum}) we may give 
an analytic condition for $S$ to be 2-free:

\begin{lemma}
\label{to be 2-free}
$S$ is 2-free in $AG(n, 2)$ if, and only if, for any triple of distinct elements 
$\alpha, \beta, \gamma \in S, \quad \alpha +\beta +\gamma \not \in S $.
\end{lemma}

\begin{theorem}
Let $S$ be a 2-free set of maximal order in $AG(n, 2)$. Then
$$
c \cdot 2^{n/2}+1 \leq |S| < \sqrt{3}\cdot 2^{n/2}
$$
where $c=1$ when $n=0\ (mod\ 2)$, and $c=2^{-1/2}$ when $n=1\ (mod\ 2)$.
\end{theorem}

\noindent
{\tt Proof} 

Suppose $S=\{s_1, s_2, \dots , s_m\}$ is a 2-free set. Then all the sums $s_1 + s_i + s_j$, where
$i \neq j \quad 2 \leq i, j \leq m$, are distinct and  $s_1 + s_i + s_j \not \in S$, 
by Lemma (\ref{to be 2-free}). Then

$$
{{m-1} \choose 2} + m = \frac{(m-1)(m-2)}{2} + m \leq 2^n, \quad \mbox{and}
$$
$$
m < \sqrt{3} \cdot 2^{n/2}
$$

To get the lower bound first suppose that $n=0 \ (mod \ 2)$ and consider $\Omega$ - the roots of unity of
degree $2^{n/2} + 1$. Easy to note that $|\Omega|=2^{n/2} + 1$ since $f(x)=x^{2^{n/2}+1}-1=0$ completely
splits in $F_{2^n}^*$, and it has no multiple roots for $f'(x)=x$. 
$\Omega$ is a cyclic subgroup of $F_{2^n}^*$; let $\omega$ be a primitive root of
unity of degree $2^{n/2} + 1$ then $\Omega = \{ 1, \omega, \omega^2, \dots, \omega^{2^{n/2}}\}$.

We claim that $\Omega$ is 2-free. If not then there are distinct $x, y, z \in \Omega$ such that
$x + y + z \in \Omega$. It follows that $(x + y + z)^{2^{n/2}+1}=1$. Then
$$
1=(x + y + z)^{2^{n/2}}(x + y + z)=(x^{2^{n/2}} + y^{2^{n/2}} + z^{2^{n/2}}) \ (x + y + z), 
$$
and since $x, y, z \in \Omega$ 
$$
x^{2^{n/2}}=x^{-1}, \quad y^{2^{n/2}}=y^{-1}, \quad z^{2^{n/2}}=z^{-1}. 
$$
So
$$
(x^{-1} + y^{-1} + z^{-1})(x+ y + z)=1
$$
Let $x=\omega^\alpha$, $y=\omega^\beta$, $z=\omega^\gamma$ where $0\leq \alpha < \beta < \gamma \leq 2^{n/2}$.
Then 
$$
\omega^{\alpha-\beta}+\omega^{\beta - \alpha}+\omega^{\alpha-\gamma}+\omega^{\gamma - \alpha}+
\omega^{\beta - \gamma}+\omega^{\gamma - \beta}=0.
$$ 
Denote $p=\beta-\alpha$, $q=\gamma-\alpha$. It is clear that $p<q$ and $0<p<q \leq 2^{n/2}$. Now we have
$$
\omega^{-p}+ \omega^p + \omega^{-q} + \omega^q + \omega^{p-q} + \omega^{q-p}=0.
$$ 
Multiplying by $\omega^{p+q}$ we get
$$
(\omega^p+ \omega^q) +\omega^{p+q} (\omega^p+ \omega^q)+(\omega^p+ \omega^q)^2=0,
$$
$$
(\omega^p+ \omega^q) (1+\omega^p) (1+ \omega^q)=0.
$$ 
This is a contradiction to our assumption $0<p<q \leq 2^{n/2}$. So $\Omega$ is 2-free.
\noindent
When $n=1 \ (mod \ 2)$ we consider $\Omega$ in $F_{2^{n-1}}$ and then take its linear embedding 
into $F_{2^n}$. 

\noindent
Summarizing, if $S$ is a 2-free set of maximal order in $AG(n, 2)$, then

$$
2^{n/2}+1 \leq |S| < \sqrt{3} \cdot 2^{n/2} \qquad \mbox{when $n=0 \ (mod \ 2)$}
$$
and
$$
2^{(n-1)/2}+1 \leq |S| < \sqrt{3} \cdot 2^{n/2} \qquad \mbox{when $n=1 \ (mod \ 2)$}.
$$ 
$\hfill{\Box}$

{\bf Remark} An interesting problem is to characterize the affine spaces for 
which the construction in Section 2 gives, in fact, a maximal linefree point set.

\end{document}